\documentclass[a4paper,11pt,oneside]{amsart}
\usepackage{amssymb, color, graphics, colordvi,euscript}
\usepackage{amscd}
\usepackage{graphicx}
\usepackage{epsfig}
\usepackage{xypic}

\newcommand{\bb}{\mathbb}

\newcommand{\rr}{\bb R}
\newcommand{\R}{\rr}

\renewcommand{\P}{\mathbb{P}}

\newcommand{\Spec}{\operatorname{Spec}}
\newcommand{\Diff}{\operatorname{Diff}}

\newcommand{\red}{\mathrm{red}}

\newcommand{\SC}{\mathcal{C}}

\newcommand{\SE}{\mathcal{E}}
\newcommand{\SF}{\mathcal{F}}

\newcommand{\SY}{\mathcal{Y}}

\newcommand{\SR}{\mathcal{R}}
\newcommand{\SX}{\mathcal{X}}

\newcommand{\be}{\mathbf{e}}

\newcommand{\aut}{\operatorname{Aut}}

\let\ra\rightarrow

\let\lra\longrightarrow

\newcommand{\id}{\mathrm{id}}

\newcommand{\autalg}{\aut}

\renewcommand{\phi}{\varphi}
\renewcommand{\tilde}{\widetilde}

\numberwithin{equation}{section}%
\newtheorem{theo}[equation]{Theorem}%
\newtheorem{thm}[equation]{Theorem}%

\newtheorem{lem}[equation]{Lemma}

\newtheorem{defin}[equation]{Definition}
\newtheorem{dfn}[equation]{Definition}
\newtheorem{question}[equation]{Question}
\theoremstyle{remark}

\newtheorem{example}[equation]{Example}

\newtheorem*{ack}{Acknowledgments}

\begin{document}

\title[Weighted blow-up singularities on real rational
surfaces]{Automorphisms of real rational 
surfaces and weighted blow-up singularities}

\author{ Johannes Huisman \and Fr\'ed\'eric Mangolte }
\address{Johannes Huisman, D\'epartement de Math\'ematiques,
Laboratoire CNRS UMR 6205, Universit\'e de Bretagne Occidentale, 6,
avenue Victor Le Gorgeu,\break CS~93837, 29238 Brest cedex 3, France.
Tel.: +33 (0)2 98 01 61 98, Fax: +33 (0)2 98 01 67 90}
\email{johannes.huisman@univ-brest.fr}
\urladdr{http://pageperso.univ-brest.fr/$\sim$huisman}
\address{Fr\'ed\'eric Mangolte, Laboratoire de Math\'ematiques,
Universit\'e de Savoie, 73376 Le Bourget du Lac Cedex, France.
Tel.: +33 (0)4 79 75 86 60, Fax: +33 (0)4 79 75 81 42}
\email{mangolte@univ-savoie.fr}
\urladdr{http://www.lama.univ-savoie.fr/$\sim$mangolte}

\thanks{The research of the second named author was partially supported by the  
ANR grant "JCLAMA" of the french "Agence Nationale de la Recherche". He benefitted also 
from the hospitality of the University of Princeton.}

\date{}

\begin{abstract}
Let $X$ be a singular real rational surface obtained from a
smooth real rational surface by performing weighted blow-ups.
 Denote by~$\autalg(X)$ the group of algebraic
 automorphisms of~$X$ into itself.
 Let~$n$ be a natural integer and let~$\be=[e_1,\ldots,e_\ell]$ be a
 partition of~$n$. Denote by~$X^\be$ the set of
 $\ell$-tuples~$(P_1,\ldots,P_\ell)$ of distinct nonsingular curvilinear
 infinitely near points of~$X$ of orders~$(e_1,\ldots,e_\ell)$. We
 show that the group~$\autalg(X)$ acts transitively on~$X^\be$.
 This statement generalizes earlier work where the case of the
 trivial partition~$\be=[1,\ldots,1]$ was treated under the
 supplementary condition that~$X$ is nonsingular.

 As an application we classify singular real rational surfaces
obtained from nonsingular surfaces by performing
weighted blow-ups.
\end{abstract}

\maketitle

\begin{quote}\small
\end{quote}

\section{Introduction}\label{sec:intro}

Let~$X$ be a nonsingular compact connected real algebraic surface,
i.e.~$X$ is a nonsingular compact connected real algebraic subset of
some~$\R^m$ of dimension~$2$. Recall that~$X$ is \emph{rational} if
the field of rational functions~$\R(X)$ of~$X$ is a purely
transcendent field extension of~$\R$ of transcendence degree~$2$.
More geometrically, $X$ is rational if there are nonempty Zariski open
subsets $U$~and $V$ of~$\R^2$ and $X$, respectively, such that there
is an isomorphism of real algebraic varieties|in the sense
of~\cite{BCR}|between $U$~and $V$.  Loosely speaking, $X$ is rational
if a nonempty Zariski open subset of~$X$ admits a rational
parametrization by a nonempty Zariski open subset of~$\R^2$ 
(The last condition imposes a priori only that $X$ is unirational but in dimension 2, 
unirationality implies rationality).  A
typical example of a rational compact real algebraic surface
is the unit sphere~$S^2$ in~$\R^3$. A rational parametrization in that
case is the inverse stereographic projection.

It has recently been shown that any rational nonsingular compact
real algebraic surface is isomorphic either to the real
algebraic torus $S^1\times S^1$, or to a real algebraic surface
obtained from the real algebraic sphere~$S^2$ by blowing up a finite
number of points~\cite{BH07,hm1}.

In the sequel, it will be convenient to identify the real algebraic
surface~$X$ with the affine scheme $\Spec\SR(X)$, where~$\SR(X)$
denotes the $\R$-algebra of all algebraic|also called
regular|functions on~$X$~\cite{BCR}.  A real-valued function~$f$
on~$X$ is \emph{algebraic} if there are real polynomials $p$~and $q$
in~$x_1,\ldots,x_m$ such that $q$ does not vanish on~$X$ and such
that~$f=p/q$ on~$X$. The algebra $\SR(X)$ is the localization of the
coordinate ring~$\R[x_1,\ldots,x_m]/I$ with respect to the
multiplicative system of all polynomials that do not vanish on~$X$,
where $I$ denotes the vanishing ideal of~$X$.  It is the subring
of~$\R(X)$ of all rational functions on~$X$ that do not have any poles
on~$X$.

Thanks to the above convention, we can define a {\em curvilinear subscheme} of~$X$ to be a closed subscheme~$P$ of~$X$ that is isomorphic
to~$\Spec\R[x]/(x^e)$, for some nonzero natural integer~$e$. We
call~$e$ the \emph{length} or the \emph{order} of~$P$. 
Let~$P$ be a curvilinear subscheme of~$X$. The reduced
scheme~$P_\red$ is an ordinary point of~$X$.  

 A curvilinear subscheme of~$X$ is also called a \emph{curvilinear infinitely
near point over $P_\red$} or an \emph{infinitesimal arc} based on~$P_\red$.
A curvilinear infinitely
near point of length~$1$ is an ordinary point of~$X$.
A curvilinear infinitely
near point of length~$2$ is a pair~$(P,L)$, where $P$ is
a point of~$X$ and $L$ is $1$-dimensional subspace of the tangent
space~$T_PX$ of~$X$ at~$P$. Equivalently, a curvilinear infinitely
near point~$P$
of~$X$ of length~$2$ is a point of the exceptional divisor of the real
algebraic surface obtained by blowing up~$X$ in an ordinary point. By
induction, a curvilinear infinitely
near point of~$X$ of length~$e$ is a
curvilinear infinitely
near point of length~$e-1$ on the exceptional divisor~$E$
of a blow-up of~$X$ (cf.~\cite[p.~171]{Mumford}).

Let~$P$ and~$Q$ be
curvilinear subschemes of~$X$. We say that $P$ and $Q$ are
\emph{distant} if the points $P_\red$ and $Q_\red$ of~$X$ are
distinct.

Let~$n$ be a natural integer and let~$\be=[e_1,\ldots,e_\ell]$ be a
partition of~$n$, where~$\ell$ is some natural integer. Denote
by~$X^\be$ the set of $\ell$-tuples~$(P_1,\ldots,P_\ell)$ of mutually distant
curvilinear subschemes~$P_1,\ldots,P_\ell$ of~$X$ of
orders~$e_1,\ldots,e_\ell$, respectively.

Recall that an \emph{algebraic
automorphism} of~$X$ is a bijective map~$f$ from~$X$ into itself
such that all coordinate functions of $f$~and $f^{-1}$ are algebraic
functions on~$X$~\cite{hm1}.  Denote by~$\autalg(X)$ the group of
algebraic automorphisms of~$X$ into itself. Equivalently,
$\autalg(X)$ is the group of~$\R$-algebra automorphisms
of~$\SR(X)$.  One has a natural action of~$\autalg(X)$ on~$X^\be$.
One of the main results of the paper is the following. 

\begin{thm}\label{thtrans}
Let~$X$ be a nonsingular rational compact real algebraic
surface. Let~$n$ be a natural integer and
let~$\be=[e_1,\ldots,e_\ell]$ be a partition of~$n$, for some
natural integer~$\ell$.  Then the group~$\autalg(X)$ acts
transitively on~$X^\be$.
\end{thm}

Roughly speaking, Theorem~\ref{thtrans} states that the
group~$\autalg(X)$ acts $\ell$-transitively on curvilinear infinitely near
points of~$X$, for any~$\ell$.  The statement generalizes earlier work
where $l$-transitivity was proved for ordinary points only, i.e., in
case of the trivial partition~$\be=[1,\ldots,1]$ (cf.~\cite{hm1}).

The statement of Theorem~\ref{thtrans} motivates the following question.

\begin{question}
 Let~$X$ be a nonsingular rational compact real algebraic
 surface. Is the subset~$\autalg(X)$ of algebraic automorphisms
 of~$X$ dense in the set~$\Diff(X)$ of all diffeomorphisms of~$X$?
 Equivalently, can any diffeomorphism of~$X$ be approximated by
 algebraic automorphisms?
\end{question}

This question is studied in the forthcoming paper \cite{km}.

The problem of approximating smooth maps between real algebraic
varieties by algebraic maps has been studied by numerous
authors~\cite{BK1,BK2,BKS,Ku,JK03,JM04,Ma06}.

It should be noted that Theorem~\ref{thtrans} does not seem to follow
from the known $n$-transitivity of~$\autalg(X)$ on ordinary points. The
difficulty is that if two $n$-tuples $P$~and $Q$ of ordinary distinct
points of~$X$ tend to two $\ell$-tuples of mutually distant curvilinear infinitely near
points with lengths $e_1,\ldots,e_\ell$, then the algebraic
diffeomorphisms mapping $P$ to~$Q$ do not necessarily have a limit
in~$\autalg(X)$.

Our proof of Theorem~\ref{thtrans} goes as follows.  First we show
that the statement of Theorem~\ref{thtrans} is valid for the real
algebraic surfaces $S^2$~and $S^1\times S^1$ by explicit construction
of algebraic automorphisms (see Theorems \ref{thsphere}~and
\ref{thtorus}). This step does use the earlier work mentioned above.
Then, we use the fore-mentioned fact that an arbitrary nonsingular
rational compact real algebraic surface is either
isomorphic to~$S^1\times S^1$, or to a real algebraic
surface obtained from~$S^2$ by blowing up a finite number of distinct
ordinary points~\cite[Theorem~4.3]{hm1}.

In order to give an application of Theorem~\ref{thtrans}, we need to
recall the following. Let~$X$ be a nonsingular rational compact
connected real algebraic surface, and let~$P$ be a curvilinear infinitely near
point of~$X$.  The \emph{blow-up} of~$X$ at~$P$ is the
blow-up~$B_P(X)$ of~$X$ at the sheaf of ideals defined by the closed
subscheme~$P$. Explicitly, if $P$ is defined by the ideal~$(x^e,y)$ on
the real affine plane~$\R^2$, then the blow-up of~$\R^2$ at~$P$ is the
real algebraic sub-variety of~$\R^2\times\P^1(\R)$ defined by the
equation~$vx^e-uy=0$, where $(u:v)$ are homogeneous coordinates on the
real projective line~$\P^1(\R)$. The blow-up~$B_P(X)$ is also called a
\emph{weighted blow-up}, for obvious reasons.  If~$e=1$, the
blow-up~$B_P(X)$ is the ordinary blow-up of~$X$ at~$P$. If~$e\geq 2$
then the blow-up~$B_P(X)$ has a singular point. A local equation of
the singularity is~$x^e=uy$ in~$\R^3$. This is often called a
singularity {\em of type~$A_{e-1}^-$} (see
e.g.~\cite[Definition~2.1]{Ko00}).

Note that weighted blow-ups recently turned out to have several
applications in real algebraic geometry
(see~\cite{KoIII,Ko00,CM08a,CM08b}).

\medskip
We apply our results, and study singular real rational surfaces
that are obtained from nonsingular ones by performing weighted blow-ups.
The latter surfaces get horns and appear as rather diabolic to us, and will be
referred to, for brevity, by the following term.
\begin{defin}
A singular real compact surface $X$ is \emph{Dantesque} if it is 
obtained from a  nonsingular real algebraic
surface~$Y$ by performing a finite number of weighted blow-ups.  
\end{defin}

The following statement implies that any rational Dantesque surface is
obtained from $S^1\times S^1$ or~$S^2$ by blowing-up a finite number
of mutually distant curvilinear subschemes of $S^1\times S^1$
or~$S^2$, respectively. Note that in particular, such a surface is connected.

\begin{thm}\label{thnonsuccessive}
Let~$X$ be a rational Dantesque surface. Then
\begin{itemize}
\item there are mutually distant curvilinear subschemes~$P_1,\ldots,P_\ell$
on~$S^1\times S^1$ such that~$X$ is isomorphic to the real algebraic
surface obtained from~$S^1\times S^1$ by blowing
up~$P_1,\ldots,P_\ell$, or
\item there are mutually distant curvilinear subschemes~$P_1,\ldots,P_\ell$
on~$S^2$ such that~$X$ is isomorphic to the real algebraic surface
obtained from~$S^2$ by blowing up~$P_1,\ldots,P_\ell$.
\end{itemize}
\end{thm}

On a singular surface, a curvilinear infinitely near point $P$ is \emph{nonsingular} 
if~$P_\red$ is a nonsingular point.

Let~$X$ be a rational Dantesque surface.  
Let~$n$ be a natural integer and
let~$\be=[e_1,\ldots,e_\ell]$ be a partition of~$n$, where~$\ell$ is
some natural integer. Denote by~$X^\be$ the set of
$\ell$-tuples~$(P_1,\ldots,P_\ell)$ of mutually distant nonsingular curvilinear 
subschemes~$P_1,\ldots,P_\ell$ of~$X$ of orders~$e_1,\ldots,e_\ell$,
respectively.

Denote again by~$\autalg(X)$ the group of algebraic automorphisms
of the possibly singular real algebraic surface~$X$.  Note that the
definition of algebraic automorphism above makes perfectly sense for
singular varieties.  Alternatively, one can define an algebraic
automorphism of~$X$ to be an automorphism of~$\Spec\SR(X)$.  Anyway,
one has a natural action of~$\autalg(X)$ on~$X^\be$.

We have the following generalization of Theorem~\ref{thtrans} above.

\begin{thm}\label{thtranssing}
Let~$X$ be a rational Dantesque surface. 
Let~$n$ be a
natural integer and let~$\be=[e_1,\ldots,e_\ell]$ be a partition
of~$n$, for some natural integer~$\ell$.  Then the
group~$\autalg(X)$ acts transitively on~$X^\be$.
\end{thm}

As an application of Theorem~\ref{thtrans} and Theorem~\ref{thnonsuccessive},
we prove the following statement.

\begin{thm}\label{thiso}
Let~$n$ be a natural integer and let~$\be=[e_1,\ldots,e_\ell]$ be a
partition of~$n$, for some natural integer~$\ell$.  Let~$X$ and $Y$
be two rational Dantesque surfaces. 
Assume that each of the surfaces
$X$~and $Y$ contains exactly one singularity of type $A_{e_i}^-$ for each
$i=1,\ldots,\ell$.  Then $X$~and $Y$ are isomorphic as real
algebraic surfaces if and only if they are homeomorphic as singular
topological surfaces.
\end{thm}

Theorem~\ref{thiso} generalizes to certain singular real rational surfaces an earlier result for nonsingular ones
(\cite[Theorem~1.2]{BH07} and~\cite[Theorem~1.5]{hm1}).
We will show by an example that the statement of Theorem~\ref{thiso}
does not hold for the slightly more general class of real rational compact
surfaces that contain singularities of
type~$A^-$ (see Example~\ref{exiso}).

\begin{ack}
We want to thank L.~Evain,  J.~Koll\'ar,  and D.~Naie for useful discussions.
\end{ack}

\section{Infinitely near points on the torus}
\label{setorus}

The object of this section is to prove Theorem~\ref{thtrans} in the
case of the real algebraic torus:

\begin{thm}\label{thtorus}
Let~$n$ be a natural integer and let~$\be=[e_1,\ldots,e_\ell]$ be a
partition of~$n$, for some natural integer~$\ell$. The
group~$\autalg(S^1\times S^1)$ acts transitively on~$(S^1\times
S^1)^\be$.
\end{thm}

The above statement is a generalization of the following statement, that
we recall for future reference.

\begin{thm}[\protect{\cite[Theorem~1.3]{BH07}}]\label{thntorus}
Let~$n$ be a natural integer.  The group $\autalg(S^1\times S^1)$
acts $n$-transitively on~$S^1\times S^1$.\qed
\end{thm}

For the proof of Theorem~\ref{thtorus}, we need several lemmas.  It
will turn out to be convenient to replace~$S^1$ by the isomorphic real
projective line~$\P^1(\R)$.

\begin{lem}\label{lephitorus}
Let~$p,q\in\R[x]$ be real polynomials in~$x$ of the same degree.
Suppose that~$q$ does not have any real roots. Define
\[
\phi\colon \P^1(\R)\times\P^1(\R)\lra\P^1(\R)\times\P^1(\R)
\]
by
\[
\phi(x,y)=\left(x,y+\frac{p}{q}\right).
\]
Then $\phi$ is an algebraic automorphism of~$\P^1(\R)\times\P^1(\R)$
into itself.
\end{lem}

\begin{proof}
It suffices to prove that~$\phi$ is an algebraic map. We
write~$\phi$ in bi-homogeneous coordinates:
\[
\phi([x_0:x_1],[y_0:y_1])=
\left([x_0:x_1],
[\overline{q}(x_0,x_1)y_0+\overline{p}(x_0,x_1)y_1: \overline{q}(x_0,x_1)y_1]
\right),
\]
where $\overline{p}$ and $\overline{q}$ are the homogenizations of
$p$~and $q$, respectively. Since~$q$ has no real zeros, the
homogeneous polynomial~$\overline{q}$ does not vanish on~$\P^1(\R)$.
Therefore, if
\begin{align*}
\overline{q}(x_0,x_1)y_0+\overline{p}(x_0,x_1)y_1&=0, \text{ and}\\
\overline{q}(x_0,x_1)y_1&=0
\end{align*}
then~$y_1=0$~and $y_0=0$. It follows that~$\phi$ is a well defined
algebraic map from~$\P^1(\R)\times\P^1(\R)$ into itself.
\end{proof}

\begin{dfn}
Let~$P$ be a curvilinear infinitely
near point of a nonsingular compact real algebraic
surface~$X$. If~$R$ is an ordinary
point of~$X$, then $P$ is said to be \emph{infinitely near} to~$R$
if~$P_\red=R$.
\end{dfn}

\begin{dfn}
Let~$P$ be a curvilinear infinitely
near point of~$S^1\times S^1$.  We say
that~$P$ is \emph{vertical} if~$P$ is tangent to a vertical
fiber~$\{x\}\times S^1$, for some~$x\in S^1$, i.e. if the
scheme-theoretic intersection~$P\cdot(\{x\}\times S^1)$ is not
reduced.
\end{dfn}

\begin{lem}\label{lenonvertical}
Let~$P_1,\ldots,P_\ell$ be mutually distant curvilinear subschemes
of~$\P^1(\R)\times \P^1(\R)$.  Then there is an algebraic
diffeomorphism~$\phi$ of~$\P^1(\R)\times\P^1(\R)$ such that
\begin{itemize}
\item $\phi(P_i)$ is infinitely near to the point~$(i,0)$
of~$\P^1(\R)\times\P^1(\R)$, and
\item $\phi(P_i)$ is not vertical,
\end{itemize}
for all~$i$.
\end{lem}

\begin{proof}
By Theorem~\ref{thntorus}, we may assume
that~$P_1,\ldots,P_\ell$ are infinitely near to the
points~$(1,0),\ldots,(\ell,0)$ of the real algebraic
torus~$\P^1(\R)\times\P^1(\R)$, respectively.

Let~$v_i=(a_i,b_i)$ be a tangent vector to~$\P^1(\R)\times\P^1(\R)$
at~$P_i$ that is tangent to~$P_i$. This means the following.
If~$P_i$ is an ordinary point then $v_i=0$. If $P_i$ is not an ordinary
point then $v_i\neq0$ and the $0$-dimensional subscheme
of~$\P^1(\R)\times\P^1(\R)$ of length~$2$ defined by~$v_i$ is
contained in the closed subscheme~$P$ of~$\P^1(\R)\times\P^1(\R)$.

Let~$p,q\in\R[y]$ be real polynomials in~$y$ of the same degree
such that
\begin{itemize}
\item $q$ does not have any real roots, 
\item $p(0)=0$, $q(0)=1$, $q'(0)=0$ and
\item $a_i+b_ip'(0)\neq0$ whenever~$v_i\neq0$.
\end{itemize}

Define
$\phi\colon \P^1(\R)\times\P^1(\R)\lra\P^1(\R)\times\P^1(\R)$ by
\[
\phi(x,y)=\left(x+\frac{p}{q},y\right).
\]
According to Lemma~\ref{lephitorus}|exchanging $x$~and $y$|the map
$\phi$ is an algebraic automorphism of~$\P^1(\R)\times\P^1(\R)$.
Since~$p(0)=0$, one has~$\phi((i,0))=(i,0)$. It follows
that~$\phi(P_i)$ is also infinitely near to~$(i,0)$.

The Jacobian of~$\phi$ at~$(i,0)$ is equal to
\[
D_{(i,0)}\phi=\begin{pmatrix}
1&\frac{p'(0)q(0)-p(0)q'(0)}{q(0)^2}\\
0&1
\end{pmatrix}=\begin{pmatrix}
1&p'(0)\\
0&1
\end{pmatrix}
\]
By construction, $(D_{(i,0)}\phi) v_i$ has first coordinate non zero
whenever~$v_i\neq0$.  Therefore, $\phi(P_i)$ is not vertical, for
all~$i$.
\end{proof}

\begin{proof}[Proof of Theorem~\ref{thtorus}] Let~$P_1,\ldots,P_\ell$
be mutually distant curvilinear subschemes of the real algebraic
torus~$\P^1(\R)\times\P^1(\R)$ of orders $e_1,\ldots,e_\ell$,
respectively.  Let~$Q_i$ be the curvilinear infinitely
near point
of~$\P^1(\R)\times\P^1(\R)$ defined by the ideal~$((x-i)^{e_i},y)$
in~$\SR(\P^1(\R)\times\P^1(\R))$. It suffices to show that there is
an algebraic automorphism~$\phi$ of~$\P^1(\R)\times\P^1(\R)$ such
that~$\phi(Q_i)=P_i$ for all~$i$.

By Lemma~\ref{lenonvertical}, we may assume that the curvilinear infinitely
near point~$P_i$ is infinitely near to~$(i,0)$ and
that~$\phi(P_i)$ is not vertical. It follows that~$P_i$ is defined
by an ideal of the form
\[
((x-i)^{e_i},y-f_i),
\]
where $f_i\in\R[x]$.

Let~$p,q\in\R[x]$ be of the same degree such that
\begin{itemize}
\item $q$ does not have any real roots, 
\item $p=f_iq$ modulo~$(x-i)^{e_i}$ for all~$i$.
\end{itemize}
Such polynomials abound by the Chinese Remainder Theorem.

By Lemma~\ref{lephitorus}, the polynomials $p$~and $q$ give rise to an
algebraic automorphism $\phi$ of~$\P^1(\R)\times\P^1(\R)$ defined by
\[
\phi(x,y)=\left(x,y+\frac{p}{q}\right).
\]

In order to show that~$\phi(Q_i)=P_i$ for all~$i$,
we compute
\[ 
(\phi^{-1})^\star((x-i)^{e_i})=(x-i)^{e_i}
\]
and
\[
(\phi^{-1})^\star(y)=y-\frac{p}{q}=y-f_i
\]
modulo~$(x-i)^{e_i}$. Indeed, $q$ is invertible modulo~$(x-i)^{e_i}$, and
$p=f_iq$ modulo $(x-i)^{e_i}$, by construction of $p$~and $q$.
It follows that~$\phi(Q_i)=P_i$.
\end{proof}

\section{Infinitely near points on the unit sphere}
\label{sesphere}

The object of this section is to prove Theorem~\ref{thtrans} in the
case of the real algebraic sphere~$S^2$:

\begin{thm}\label{thsphere}
Let~$n$ be a natural integer and let~$\be=[e_1,\ldots,e_\ell]$ be a
partition of~$n$, for some natural integer~$\ell$. The
group~$\autalg(S^2)$ acts transitively on~$(S^2)^\be$.
\end{thm}

The above statement is a generalization of the following statement, that
we recall for future reference.

\begin{theo}[\protect{\cite[Theorem~2.3]{hm1}}]\label{thnsphere}
Let $n$ be a natural integer. The group $\autalg(S^2)$ acts
$n$-transitively on $S^2$.\qed
\end{theo}

For the proof of Theorem~\ref{thsphere}, we need several lemmas. 

\begin{lem}[\protect{\cite[Lemma~2.1]{hm1}}]\label{lephisphere}
Let~$p,q,r\in\R[x]$ be such that
\begin{itemize}
\item $r$ does not have any roots in the interval~$[-1,1]$, and
\item $p^2+q^2=r^2$.
\end{itemize}
Define $\phi\colon S^2\lra S^2$
by
\[
\phi(x,y,z)=\left(x,\frac{yp-zq}{r},\frac{yq+zp}{r}\right).
\]
Then $\phi$ is an algebraic automorphism of~$S^2$.\qed
\end{lem}

\begin{dfn}
Let~$P$ be infinitely near to a point of the equator~$\{z=0\}$
of~$S^2$.  We say that~$P$ is \emph{vertical} if~$P$ is tangent to the
great circle of~$S^2$ passing through the North pole.
\end{dfn}

As for the torus above, we need some standard points on~$S^2$.  Let
\[
R_i=(x_i,y_i,z_i)=\left(\cos({\textstyle\frac{i\pi}{2\ell+1}}),
\sin({\textstyle\frac{i\pi}{2\ell+1}}),0\right)
\]
for $i=1,\ldots,\ell$. Note that~$x_i\neq0$ and $y_i\neq0$ for all~$i$.

\begin{lem}\label{lenonhorizontal}
Let~$P_1,\ldots,P_\ell$ be mutually distant curvilinear subschemes
of~$S^2$.  Then there is an algebraic automorphism~$\phi$
of~$S^2$ such that
\begin{itemize}
\item $\phi(P_i)$ is infinitely near to~$R_i$, and
\item $\phi(P_i)$ is not vertical, 
\end{itemize}
for all~$i$.
\end{lem}

\begin{proof}
By Theorem~\ref{thnsphere}, we may assume that~$P_i$ is
infinitely near to the point~$R_i$, for all~$i$.
Let~$v_i=(a_i,b_i,c_i)$ be a tangent vector to~$S^2$ at~$P_i$ that
is tangent to~$P_i$.  Let~$p,q,r\in\R[x]$ be such that
\begin{itemize}
\item $r$ does not vanish on~$[-1,1]$,
\item $p^2+q^2=r^2$,
\item $p(0)=1$, $q(0)=0$, $r(0)=1$, $p'(0)=0$, $r'(0)=0$, and
\item $a_i-c_iy_iq'(0)\neq0$ or $b_i+c_ix_iq'(0)\neq0$,
 whenever~$v_i\neq0$.
\end{itemize}
Such polynomials abound.  Take, for example,
$$ p(z)=(1+z^2)^2-(\lambda z)^2,\quad q(z)=2(1+z^2)\lambda z, \quad
r(z)=(1+z^2)^2+(\lambda z)^2,
$$
where~$\lambda$ is any real number such that~$a_i-2\lambda y_ic_i\neq0$
or $b_i+2\lambda x_ic_i\neq0$ whenever~$v_i\neq0$.

Define $\phi\colon S^2\lra S^2$ by
\[
\phi(x,y,z)=\left(\frac{xp(z)-yq(z)}{r(z)},\frac{xq(z)+yp(z)}{r(z)},z\right).
\]
According to Lemma~\ref{lephisphere}|permuting the roles
of~$x,y,z$|the map $\phi$ is an algebraic automorphism of~$S^2$.
Since $p(0)=1$, $q(0)=0$~and $r(0)=1$, the curvilinear infinitely
near point
$\phi(P_i)$ is again infinitely near to~$R_i$, for all~$i$.

The Jacobian of~$\phi$ at~$R_i$ is equal to
\begin{multline*}
D_{R_i}\phi=\begin{pmatrix}
\frac{p(0)}{r(0)} & \frac{-q(0)}{r(0)} & \frac{x_ip'(0)r(0)-y_iq'(0)r(0)-x_ip(0)r'(0)+y_iq(0)r'(0)}{r(0)^2}\\
\frac{q(0)}{r(0)} & \frac{p(0)}{r(0)} & \frac{x_iq'(0)r(0)+y_ip'(0)r(0)-x_iq(0)r'(0)-y_ip(0)r'(0)}{r(0)^2}\\
0&0&1 
\end{pmatrix}=\\
\begin{pmatrix}
1&0&-y_iq'(0)\\
0&1&x_iq'(0)\\
0&0&1 
\end{pmatrix}
\end{multline*}
By construction, $(D_{R_i}\phi) v_i$ has first or second coordinate
non zero whenever~$v_i\neq0$.  Therefore, $\phi(P_i)$ is not
vertical, for all~$i$.
\end{proof}

\begin{lem}\label{lemod}
Let~$e$ be a nonzero natural integer, and
let~$i\in\{1,\ldots,\ell\}$. Let~$f,g,h\in\R[x]/(x-x_i)^e$ be such
that
\begin{equation}\label{eqmod}
x^2+f^2=1,\quad\text{and}\quad
x^2+g^2+h^2=1
\end{equation}
in~$\R[x]/(x-x_i)^e$. Assume, moreover, that~$f(x_i)=g(x_i)=y_i$.
Then there is $a\in\R[x]/(x-x_i)^e$ such that
\begin{equation}\label{eqmod1}
(1-a^2)f=(1+a^2)g,\quad\text{and}\quad
2af=(1+a^2)h
\end{equation}
in~$\R[x]/(x-x_i)^e$. Moreover, there is such an element~$a$ such
that~$1+a^2$ is invertible in~$\R[x]/(x-x_i)^e$.
\end{lem}

\begin{proof}
If~$h=0$ then $f=g$, and one can take~$a=0$. Therefore, we may
assume that~$h\neq0$.  Let~$d$ be the valuation of~$h$, i.e.,
$h=(x-x_i)^dh'$, where $h'\in\R[x]/(x-x_i)^e$ is invertible.
Since~$f(x_i)=g(x_i)$, one has~$h(x_i)=0$, i.e., $h$ is not
invertible in~$\R[x]/(x-x_i)^e$ and $d\neq0$. By Hensel's Lemma,
there are lifts~$\overline{f},\overline{g},\overline{h}$
in~$\R[x]/(x-x_i)^{e+2d}$ of~$f,g,h$, respectively, satisfying the
equations~(\ref{eqmod}) in the ring~$\R[x]/(x-x_i)^{e+2d}$.  Note
that~$\overline{f}+\overline{g}$ is invertible
in~$\R[x]/(x-x_i)^{e+2d}$, and that~$\overline{h}$ has
valuation~$d$.

In order to simplify notation, we denote again by~$f,g,h$ the
elements $\overline{f},\overline{g},\overline{h}$, respectively.
Let $k\in\R[x]/(x-x_i)^{e+2d}$ be the inverse of~$f+g$.  Let~$a=hk$.
We verify that equations~(\ref{eqmod1}) hold and that $1+a^2$ is
invertible in~$\R[x]/(x-x_i)^e$.

The element~$1+a^2$ is clearly invertible in~$\R[x]/(x-x_i)^{e+2d}$
since~$h$ is not invertible.

Since
$$ 
(f-g)(f+g)=f^2-g^2=(1-x^2)-(1-x^2-h^2)=h^2,
$$
one has
$$
f-g=h^2k=h^2k^2(f+g)=a^2(f+g).
$$
It follows that
$$
(1-a^2)f=(1+a^2)g
$$
in $\R[x]/(x-x_i)^{e+2d}$, and therefore also in~$\R[x]/(x-x_i)^e$.

In order to prove that the other equation of~(\ref{eqmod1}) holds as well, 
observe that
\begin{multline*}
(f-g)^2h-2f(f-g)h+h^3=
(f-g)h(f-g-2f)+h^3=\\
-(f-g)h(f+g)+h^3=
-(f^2-g^2)h+h^3=0
\end{multline*}
by what we have seen above. Substituting~$f-g=ah$, one obtains
$$
0=a^2h^3-2afh^2+h^3=h^2(a^2h-2af+h)
$$
in~$\R[x]/(x-x_i)^{e+2d}$. Since the valuation of~$h$ is equal to~$d$,
one deduces that $a^2h-2af+h=0$ in~$\R[x]/(x-x_i)^e$. Hence,
$2af=(1+a^2)h$, as was to be proved.\end{proof}

\begin{proof}[Proof of Theorem~\ref{thsphere}] Let~$P_1,\ldots,P_\ell$
be mutually distant curvilinear subschemes  of~$S^2$ of orders $e_1,\ldots,e_\ell$,
respectively.  Let~$Q_i$ be the
curvilinear infinitely
near point of~$S^2$ defined by the
ideal
$$
((x-x_i)^{e_i},y-f_i,z)
$$ 
in~$\R[x,y,z]$, where $f_i$ is the Taylor polynomial in~$x-x_i$ of
$\sqrt{1-x^2}$ at~$x_i$ of order~$e_i-1$. Note that~$Q_i$ is
infinitely near to~$R_i$ for all~$i$.
We show that there is an algebraic automorphism~$\phi$ of~$S^2$
such that~$\phi(Q_i)=P_i$ for all~$i$.

By Lemma~\ref{lenonhorizontal}, we may assume
that~$P_1,\ldots,P_\ell$ are infinitely near to the points
$R_1,\ldots,R_\ell$ of~$S^2$, respectively, and that they are not
vertical.  It follows that~$P_i$ is defined by an ideal of the form
\[
((x-x_i)^{e_i},y-g_i,z-h_i)
\]
where $g_i,h_i\in\R[x]$ are of degree $<e_i$. Moreover, since~$P_i$
is a curvilinear infinitely
near point of~$S^2$, we have
$$
x^2+g_i^2+h_i^2=1\pmod{(x-x_i)^{e_i}}.
$$

By Lemma~\ref{lemod}, there is~$a_i\in\R[x]/(x-x_i)^{e_i}$ such that
\begin{align*}
(1-a_i^2)f_i&=(1+a_i^2)g_i,\quad\text{and}\\
2a_if_i&=(1+a_i^2)h_i
\end{align*}
in~$\R[x]/(x-x_i)^{e_i}$, and, moreover, $1+a_i^2$ is invertible.

By the Chinese Remainder Theorem, there is a polynomial~$a\in\R[x]$ such that 
$a=a_i  \pmod{(x-x_i)^{e_i}}$, for all~$i$. Then
\begin{align*}
(1-a^2)f_i&=(1+a^2)g_i\pmod{(x-x_i)^{e_i}}\\
2af_i&=(1+a^2)h_i\pmod{(x-x_i)^{e_i}},
\end{align*}
and $1+a^2$ is invertible in~$\R[x]/(x-x_i)^{e_i}$, for all~$i$.

Put
$$
p=1-a^2,\quad q=2a,\quad r=1+a^2.
$$  
Then
\begin{align*}
pf_i&=rg_i\pmod{(x-x_i)^{e_i}}\\
qf_i&=rh_i\pmod{(x-x_i)^{e_i}},
\end{align*}
and $r$ is invertible in~$\R[x]/(x-x_i)^{e_i}$, for all~$i$. Moreover,
\begin{itemize}
\item  $r$ does not have any roots in the interval~$[-1,1]$, and 
\item  $p^2+q^2=r^2$.
\end{itemize}
By Lemma~\ref{lephisphere}, the polynomials $p,q,r$ give rise to an
algebraic automorphism~$\phi$ of~$S^2$ defined by
\[
\phi(x,y,z)=\left(x,\frac{yp-zq}{r},\frac{yq+zp}{r}\right).
\]

In order to show that~$\phi(Q_i)=P_i$ for all~$i$, we compute
\begin{align*}
(\phi^{-1})^\star((x-x_i)^{e_i})&=(x-x_i)^{e_i}\\
u_i=(\phi^{-1})^\star(y-f_i)&=\frac{yp+zq}{r}-f_i\\
v_i=( \phi^{-1})^\star(z)&=\frac{-yq+zp}{r},
\end{align*}
so that~$\phi(Q_i)$ is the curvilinear infinitely
near point of~$S^2$ defined by
the ideal~$((x-x_i)^{e_i},u_i,v_i)$. We have
$$
\frac{p}{r}u_i-\frac{q}{r}v_i=y-\frac{p}{r}f_i=y-g_i\pmod{(x-x_i)^{e_i}}
$$
and
$$
\frac{q}{r}u_i+\frac{p}{r}v_i=z-\frac{q}{r}f_i=z-h_i\pmod{(x-x_i)^{e_i}}.
$$
It follows that~$\phi(Q_i)=P_i$.
\end{proof}

\section{algebraic automorphisms of nonsingular rational surfaces}

The proof of Theorem~\ref{thtrans} is now similar to the proof
of~\cite[Theorem~1.4]{hm1}. We include it for convenience of the reader.

\begin{proof}[Proof of Theorem~\ref{thtrans}]
Let $X$ be a nonsingular real rational surface, and let $(P_1,\dots,P_\ell)$ and $(Q_1,\dots,Q_\ell)$
be two $\ell$-tuples in $X^\be$.  As mentioned before,
$X$ is isomorphic to~$S^1\times S^1$ or to the
blow-up of~$S^2$ at a finite number of distinct
points~$R_1,\ldots,R_m$.  If $X$ is isomorphic
to~$S^1\times S^1$ then~$\autalg(X)$ acts transitively on~$X^\be$
by Theorem~\ref{thtorus}.  Therefore, we may assume that~$X$ is
isomorphic to the blow-up~$B_{R_1,\dots,R_m}(S^2)$
of~$S^2$ at~$R_1,\ldots,R_m$.  Moreover, we may assume that the points
$(P_1)_\red,\ldots,(P_n)_\red,(Q_1)_\red,\ldots,(Q_n)_\red$ do not belong to any of the
exceptional divisors, by ~\cite[Theorem~3.1]{hm1}. Thus we can consider the $P_j,Q_j$ as curvilinear subschemes of $S^2$. It
follows that $(R_1,\ldots,R_m,P_1,\ldots,P_\ell)$ and
$(R_1,\ldots,R_m,Q_1,\ldots,Q_\ell)$ are two $(m+\ell)$-tuples in~$(S^2)^{\mathbf f}$, where $\mathbf f = [1,\dots,1,e_1,\ldots,e_\ell]$.

By Theorem~\ref{thsphere}, there is an automorphism
$\psi$ of~$S^2$ such that~$\psi(R_i)=R_i$, for all~$i$,
and~$\psi(P_j)=Q_j$, for all~$j$.  The induced automorphism~$\phi$
of~$X$ has the property that~$\phi(P_j)=Q_j$, for all~$j$.

\end{proof}

\section{Rational surfaces with $A^-$ singularities}
\label{seaminus}

The object of this section is to prove Theorem~\ref{thnonsuccessive}
that asserts that a real rational weighted blow-up surface is isomorphic 
to a real algebraic surface obtained from
$S^2$ or~$S^1\times S^1$ by blowing up a finite number of mutually distant curvilinear subschemes.

The following is a particular case of~\cite[Theorem~3.1]{hm1}.

\begin{lem}\label{lekleinbottle}
Let~$X$ be a nonsingular rational real algebraic
Klein bottle Let~$S$ be a finite subset of~$X$. Then there is an
algebraic map~$f\colon X\ra S^2$ such that
\begin{enumerate}
\item $f$ is the blow-up of~$S^2$ at $2$ distinct real
points~$Q_1,Q_2$, and
\item $Q_i\not\in f(S)$, for~$i=1,2$.
\end{enumerate}
\end{lem}

\begin{lem}\label{leblow-uptorus}
Let~$P$ be a curvilinear infinitely
near point of~$S^1\times S^1$, and let~$C$
be a real algebraic curve in $S^1\times S^1$ 
such that there is a nonsingular projective complexification~$\SX$
of~$S^1\times S^1$ having the following properties:
\begin{enumerate}
\item the Zariski closure~$\SC$ of~$C$ in~$\SX$ is nonsingular and rational,
\item the self-intersection of~$\SC$ in~$\SX$ is even and non-negative,
\item $P_\red\in\SC$, and
\item $\SC$ is not tangent to~$P$, i.e., the scheme-theoretic
intersection~$P\cdot\SC$ is of length~$1$.
\end{enumerate}
Then, there is an algebraic map
$$
f\colon B_P(S^1\times S^1)\ra Z
$$
that is a blow-up at a curvilinear infinitely
near point~$Q$ whose exceptional
curve $f^{-1}(Q_\red)$ is equal to the strict transform of~$C$
in~$B_P(S^1\times S^1)$, where~$Z$ is either the real algebraic
torus~$S^1\times S^1$, or the rational real algebraic Klein bottle~$K$.
\end{lem}

\begin{proof}
Let~$Y$ be the blow-up of~$S^1\times S^1$ at~$P$.
Let~$\beta\colon\SY\ra\SX$ be the blow-up of~$\SX$ at~$P$. It is
clear that~$\SY$ is a nonsingular projective complexification
of~$Y$. 

Let~$m+1$ be the length of the curvilinear infinitely
near point~$P$,
where~$m\geq0$.  Let~$\rho\colon\tilde\SY\ra\SY$ be the minimal
resolution of~$\SY$. If~$P$ is a point of length~$1$,
then~$\rho=\id$, of course.  The morphism~$\beta\circ\rho$ is
a repeated blow-up of~$\SX$. More precisely, there is a sequence of
morphisms of algebraic varieties over~$\R$
$$
\xymatrix{
 \tilde\SY=\SX_{m+1}\ar[r]^{f_{m}}&\SX_{m}\ar[r]^{f_{m-1}}&
 \cdots\ar[r]^{f_0}&\SX_0=\SX\, , }
$$
with the following properties. Each morphism~$f_i$ is an ordinary
blow-up of~$\SX_i$ at a nonsingular ordinary real point~$P_i$
of~$\SX_i$, for all~$i$. One has $P_0=P_\red$, i.e., $f_0$ is the
blow-up of~$\SX$ at~$P_\red$. Moreover, denoting by~$\SE_i$ the
exceptional curve of~$f_i$ in~$\SX_{i+1}$, the center of
blow-up~$P_{i+1}$ belongs to~$\SE_i$ but does not belong to the
strict transform of any of the curves~$\SE_j$ in~$\SX_{i+1}$ for
all~$j<i$.

Denote again by~$\SE_i$ and~$\SC$ the strict transforms
of~$\SE_i$~and $\SC$ in~$\tilde\SY$, respectively. The
curves~$\SE_0,\ldots,\SE_{m-1}$ have self-intersection~$-2$, the
curve~$\SE_m$ has self-intersection~$-1$.  Since~$\SC$ is not
tangent to~$P$, the curve~$\SC$ in~$\tilde\SY$ has odd
self-intersection~$\geq-1$. The curves~$\SC,\SE_0,\ldots,\SE_m$ form
a chain of curves over~$\R$ in~$\SY$, intersecting in real points
only (See Figure~\ref{fichain}).
\begin{figure}
 \centering

 \epsfbox{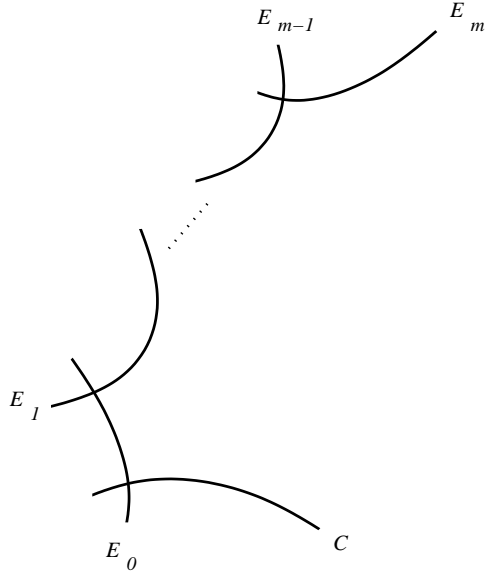} \caption{The chain of exceptional curves
   in~$\tilde{\SY}$.}
 \label{fichain}
\end{figure}
The morphism~$\rho\colon\tilde\SY\ra\SY$ is the contraction of the
curves~$\SE_0,\ldots,\SE_{m-1}$. The morphism~$\beta\colon\SY\ra\SX$
is the contraction of~$\rho(\SE_m)$, i.e.,
$\beta^{-1}(P_\red)=\rho(\SE_m)$.

Let~$k$ be the self-intersection of~$\SC$ in~$\tilde\SY$.
Since~$k\geq-1$ and $k\equiv-1\pmod2$, the integer~$k+1$ is even and
non-negative. Let~$R_1,\ldots,R_{k+1}$ be pairwise complex
conjugate of~$\SC$. Denote by~$\tilde\SY'$ the blow-up
of~$\tilde{\SY}$ in $R_1,\ldots,R_{k+1}$. The algebraic
variety~$\tilde\SY'$ is again defined over~$\R$. The strict
transform of~$\SC$ in~$\tilde\SY'$ is a nonsingular rational curve
of self-intersection~$-1$. Denote again by~$\SE_i$ the strict
transform of~$\SE_i$ in~$\tilde\SY'$. The self-intersection
of~$\SE_i$ is equal to~$-2$, if~$i\neq m$, the self-intersection
of~$\SE_m$ is equal to~$-1$.

Let~$\SY'$ be the algebraic surface defined over~$\R$ obtained
from~$\tilde\SY'$ by contracting the union of the
curves~$\SE_0\ldots,\SE_{m-1}$ to a point, and
let~$\rho'\colon\tilde\SY'\ra\SY'$ be the contracting morphism.
Let~$\SX'$ be the algebraic surface defined over~$\R$ obtained
from~$\SY'$ by contracting~$\rho'(\SC)$ to a point, and
let~$\beta'\colon\SY'\ra\SX'$ be the contracting morphism.
Since~$\beta'\circ\rho'$ is a repeated blow-down of $-1$-curves, the
algebraic surface~$\SX'$ is nonsingular.  Moreover, the
morphism~$\beta'$ is a blow-up of~$\SX'$ at a nonsingular curvilinear infinitely
near point~$Q$ of~$\SX'$. Denote again by~$\SC$ the
curve~$\rho'(\SC)$ in~$\SY'$. The curve~$\SC$ in~$\SY'$ is the
exceptional curve of~$\beta'$.

Now take the associated real algebraic varieties, denoted by the
corresponding roman characters. Since the
points~$R_1,\ldots,R_{k+1}$ are non real, one
has~$\tilde\SY'(\R)=\tilde\SY(\R)$, i.e., $\tilde Y'=\tilde Y$, the
minimal resolution of~$Y$. It follows that~$Y'=Y$, and that the
induced algebraic map~$b\colon Y'\ra X'$ is the blow-up of the
curvilinear infinitely
near point~$Q$ of the nonsingular compact connected real
algebraic surface~$X'$. The exceptional curve of~$b$ is equal to the
strict transform of~$C$ in~$Y$.

The only thing that is left to prove is the fact that the real
algebraic surface~$X'$ is isomorphic to~$S^1\times S^1$ or to
the rational real algebraic Klein bottle~$K$. In order to establish this,
observe that~$\tilde Y$, as an $(m+1)$-fold blow-up of~$S^1\times
S^1$, is homeomorphic to the connected sum of $S^1\times S^1$ and
$m+1$ copies of~$\P^2(\R)$.  Since~$Y'$ also is homeomorphic to the
connected sum of~$X'$ and $m+1$ copies of~$\P^2(\R)$, it follows
that $X'$ is homeomorphic to a torus or a Klein bottle.
By~\cite[Theorem~1.3]{Ma06}, or \cite[Theorem~1.2]{BH07}, or \cite[Theorem~1.5]{hm1}, $X'$ is
isomorphic to~$S^1\times S^1$ or the rational real
algebraic Klein bottle~$K$.
\end{proof}

A similar, but easier, argument applies and proves the following lemma.

\begin{lem}\label{leblow-upsphere}
Let~$P$ be a curvilinear infinitely
near point of~$S^2$, and let~$C$
be a real algebraic curve in $S^2$ 
such that there is a nonsingular projective complexification~$\SX$
of~$S^2$ having the following properties:
\begin{enumerate}
\item the Zariski closure~$\SC$ of~$C$ in~$\SX$ is nonsingular and rational,
\item the self-intersection of~$\SC$ in~$\SX$ is even and non-negative,
\item $P_\red\in\SC$, and
\item $\SC$ is not tangent to~$P$, i.e., the scheme-theoretic
intersection~$P\cdot\SC$ is of length~$1$.
\end{enumerate}
Then, there is an algebraic map
$$
f\colon B_P(S^2)\ra S^2
$$
that is the blow-up of~$S^2$ at a
curvilinear infinitely
near point~$Q$. Moreover, the exceptional
curve~$f^{-1}(Q)$ is equal to the strict transform of~$C$
in~$B_P(S^2)$.\qed
\end{lem}

\begin{proof}[Proof of Theorem~\ref{thnonsuccessive}]
There is a nonsingular real rational compact
surface~$Y$ such that~$X$ is isomorphic to the real algebraic
surface obtained from~$Y$ by repeatedly blowing up a nonsingular curvilinear infinitely
near point. Since~$Y$ is a nonsingular rational compact
real algebraic surface, $Y$ is obtained either from~$S^2$ or
from~$S^1\times S^1$, by repeatedly blowing up an ordinary point
(cf.~\cite[Theorem~3.1]{BH07} or \cite[Theorem~4.1]{hm1}). Hence,
there is a sequence of algebraic maps
$$
\xymatrix{
X=X_n\ar[r]^{f_n}&X_{n-1}\ar[r]^{f_{n-1}}&\cdots\ar[r]^{f_1}&X_0=Z\, ,
}
$$
where~$Z=S^2$ or~$Z=S^1\times S^1$, and each map~$f_i$ is a blow-up
at a nonsingular curvilinear infinitely
near point~$Q_i$ of~$X_{i-1}$, possibly
of length~$1$, for $i=1,\ldots,n$. Denote by~$E_i$ the exceptional
divisor~$f_i^{-1}((Q_i)_\red)$ of~$f_i$ in~$X_i$.

Let~$\SF$ be the set of the curvilinear infinitely
near points~$Q_i$.  Define a
partial ordering on~$\SF$ by $Q_i\leq Q_j$ if the composition
$f_i\circ\cdots\circ f_{j-1}$ maps~$(Q_j)_\red$ to~$(Q_i)_\red$.  It is
clear that~$\SF$ is a forest, i.e., a disjoint union of trees.

Let~$s$ be the number of edges in the forest~$\SF$.
We show the statement of the theorem by induction on~$s$. The
statement is clear for~$s=0$. Suppose, therefore, that~$s\neq0$. We
may assume that~$Q_1$ is the root of a tree of~$\SF$ of nonzero
height.

Let~$C$ be a real algebraic curve in~$Z$ satisfying the conditions
of Lemma~\ref{leblow-uptorus} if~$Z=S^1\times S^1$, and of
Lemma~\ref{leblow-upsphere} if~$Z=S^2$, with $P=(Q_1)_\red$. Such
curves abound: one can take a bi-degree~$(1,1)$ in~$S^1\times S^1$,
or a Euclidean circle in~$S^2$, respectively. Moreover, we may
assume that the strict transform of~$C$ in~$X_i$ does not
contain~$(Q_{i+1})_\red$, for all~$i\geq1$. Applying
Lemma~\ref{leblow-uptorus} and Lemma~\ref{leblow-upsphere},
respectively, one obtains a sequence
$$
\xymatrix{
X=X_n\ar[r]^{f_n}&X_{n-1}\ar[r]^{f_{n-1}}&\cdots\ar[r]^{f_2}&
X_1\ar[r]^{f_1'}&X_0'=Z'\, ,
}
$$
where $f_1'$ contracts the strict transform of~$C$ in~$X_1$ to a
point~$Q_1'$.  The real algebraic surface~$Z'$ is either the real
algebraic sphere~$S^2$, or the real algebraic torus~$S^1\times S^1$,
or the rational real algebraic Klein bottle~$K$. By construction, the
number of edges in the forest~$\SF'$ associated to the latter
sequence of blow-ups is equal to~$s-1$. Therefore, if~$Z'=S^2$ or
$Z'=S^1\times S^1$, we are done.  If~$Z$ is the real algebraic
Klein bottle~$K$, then, according to Lemma~\ref{lekleinbottle}, there is
a sequence of blow-ups
$$
\xymatrix{
Z'=X_0'\ar[r]^{f_0}&X_{-1}'\ar[r]^{f_{-1}}&X_{-2}=S^2\, 
}
$$
at ordinary points such that the images of the centers
$Q_1',Q_2,\ldots,Q_n$ in~$X_{-1}$ and~$X_{-2}$ are distinct from the
centers of the blow-ups $f_0$ and $f_{-1}$. We conclude also in this
case by the induction hypothesis.
\end{proof}

A close inspection of the above proof reveals that the following
slightly more technical statement holds.

\begin{thm}\label{thnonsuccessiveprime}
Let~$X$ be a rational Dantesque surface, 
and let~$S\subseteq X$ be a finite
subset of nonsingular points of~$X$. Then there is an algebraic map
$f\colon X\ra S^2$ or $f\colon X\ra S^1\times S^1$ with the
following properties:
\begin{enumerate}
\item there are mutually distant curvilinear subschemes~$P_1,\ldots,P_\ell$
on~$S^2$ or~$S^1\times S^1$, respectively, such that~$f$ is the
blow-up at~$P_1,\ldots,P_\ell$, and
\item $(P_i)_\red\not\in f(S)$, for all~$i$.\qed
\end{enumerate}
\end{thm}

\section{Infinitely near points on a singular rational surface}
\label{setrans}

The object of this section is to prove Theorem~\ref{thtranssing}.

\begin{proof}[Proof of Theorem~\ref{thtranssing}]
Let~$(P_1,\ldots,P_\ell)$ and $(Q_1,\ldots,Q_\ell)$ be two elements
of~$X^\be$. We prove that there is an algebraic automorphism~$\phi$
of~$X$ such that~$\phi(P_i)=Q_i$, as curvilinear infinitely
near points.

Let~$S$ be the set of ordinary points
$(P_1)_\red,\ldots,(P_\ell)_\red,(Q_1)_\red,\ldots,(Q_\ell)_\red$
of~$X$. Since~$S$ is a finite set of nonsingular points of~$X$,
there is, by Theorem~\ref{thnonsuccessiveprime}, an algebraic map
$f\colon X\ra S^2$ or $f\colon X\ra S^1\times S^1$ with the
following properties:
\begin{itemize}
\item there are mutually distant curvilinear subschemes~$R_1,\ldots,R_m$ on~$S^2$ such that~$f$ is the blow-up
at~$R_1,\ldots,R_m$, and
\item $(R_i)_\red\not\in f(S)$, for all~$i$.
\end{itemize}
Since~$f$ is an isomorphism at a neighborhood of~$S$, the
image~$f(P_i)$ is a curvilinear infinitely
near point of~$S^2$ of length~$e_i$,
and the same holds for~$f(Q_i)$, for all~$i$.

By Theorems \ref{thtorus}~and \ref{thsphere}, there is an algebraic automorphism~$\psi$ of~$S^2$ or~$S^1\times S^1$, respectively, such
that $\psi(P_i)=Q_i$ for $i=1,\ldots,\ell$, and $\psi(R_i)=R_i$ for
$i=1,\ldots,m$.  Then, $\psi$ induces an algebraic automorphism~$\phi$ of~$X$ with the required property.
\end{proof}

\section{Isomorphic rational real algebraic surfaces}

\begin{proof}[Proof of Theorem~\ref{thiso}] Let $X$~ and $Y$ be rational Dantesque surfaces, such that each of the surfaces $X$~and $Y$
contains exactly one singularity of type~$A_{e_i}^-$, for each
$i=1,\ldots,\ell$.

If~$X$ and $Y$ are isomorphic, then, of course, the singular
topological surfaces $X$~and $Y$ are homeomorphic.

Conversely, suppose that $X$~and $Y$ are homeomorphic.  By
Theorem~\ref{thnonsuccessive}, there are nonsingular real rational
surfaces $X'$~and $Y'$, and $\ell$-tuples \break
$(P_1,\ldots,P_\ell)\in (X')^\be$ and $(Q_1,\ldots,Q_\ell)\in
(Y')^\be$ such that $X$ is the blow-up of~$X'$ at $P_1,\ldots,P_\ell$
and $Y$ is the blow-up of~$Y'$ at~$Q_1,\ldots,Q_\ell$. Since~$X$ and
$Y$ are homeomorphic, $X'$~and $Y'$ are homeomorphic. It follows that
$X'$~and $Y'$ are isomorphic. By
Theorem~\ref{thtrans}, there is an isomorphism~$\phi\colon X'\ra Y'$ such that~$\phi(P_i)=Q_i$
for~$i=1,\ldots,\ell$. The isomorphism~$\phi$ induces an
isomorphism between $X$~and $Y$.
\end{proof}

The following example shows that the statement of Theorem~\ref{thiso}
does not hold for the slightly more general class of rational compact
connected real algebraic surfaces that contain singularities of
type~$A^-$.

\begin{example}\label{exiso}
Let~$X$ be the real algebraic surface obtained from the real
algebraic torus $S^1\times S^1$ by contracting a
fiber~$S^1\times\{\star\}$ to a point. Then~$X$ is a rational
compact connected real algebraic surface containing only one singular
point. Its singularity is of type~$A_1^-$. 

Let~$P$ be a point of~$\P^2(\R)$. The real algebraic surface~$K$
obtained from $\P^2(\R)$ by blowing up~$P$ is a real algebraic Klein
bottle.  Let~$Y$ be the real algebraic surface obtained from the
real algebraic Klein bottle~$K$ by contracting to a point the strict
transform of a real projective line in~$\P^2(\R)$ that passes
through~$P$. Then~$Y$ is a rational compact connected real algebraic
surface containing only one singular point. Its singularity is of
type~$A_1^-$.

It is clear that~$X$ and $Y$ are homeomorphic singular surfaces.
Indeed, they are both rational real algebraic models of the
once-pinched torus (Figure~\ref{fiopt}).
\begin{figure}
 \centering
 \epsfbox{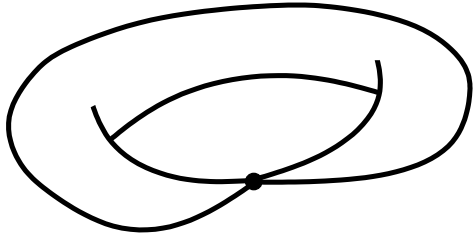} \caption{The once-pinched torus.}
 \label{fiopt}
\end{figure}
However, they are non isomorphic as real algebraic surfaces. Indeed,
if they were isomorphic, their minimal resolutions $S^1\times S^1$
and $K$ were isomorphic, which is absurd.

Note that the real rational surface~$Y$ is Dantesque, whereas $X$ is not.
\end{example}



\begin{thebibliography}{999999}

\bibitem[BH07]{BH07} 
I.~Biswas, J.~Huisman, Rational
   real algebraic models of topological surfaces, Doc.
   Math.~\textbf{12} (2007), 549--567

\bibitem[BCR98]{BCR} J.~Bochnak, M.~Coste, M.-F.~Roy,
 Real algebraic geometry, Ergeb. Math. Grenzgeb.  (3),
 vol.~36, Springer Verlag, 1998

\bibitem[BK87a]{BK1} J.~Bochnak, W.~Kucharz, Algebraic
   approximation of mappings into spheres, Michigan
   Math. J.~\textbf{34} (1987), 119--125

\bibitem[BK87b]{BK2} J.~Bochnak, W.~Kucharz,
 Realization of homotopy classes by algebraic
   mappings, J. Reine Angew. Math.~\textbf{377}
 (1987), 159--169

\bibitem[BKS97]{BKS} J.~Bochnak, W.~Kucharz, R.~Silhol,
 Morphisms, line bundles and moduli spaces in real
   algebraic geometry, Pub. Math.
   I.H.E.S.~\textbf{86} (1997), 5--65

\bibitem[CM08a]{CM08a} F.~Catanese, F.~Mangolte, Real
   singular Del Pezzo surfaces and threefolds fibred by rational
   curves, I, Michigan Math. J.~\textbf{56} (2008), 357--373

\bibitem[CM08b]{CM08b} F.~Catanese, F.~Mangolte, Real
   singular Del Pezzo surfaces and threefolds fibred by rational
   curves, II, Ann. Sci. E.~N.~S. (to appear), \texttt{arXiv:0803.2074 [math.AG]}

\bibitem[HM08]{hm1} J.~Huisman, F.~Mangolte,
The group of
  automorphisms of a real rational surface is
$n$-transitive, (submitted),
\texttt{arXiv:0708.3992 [math.AG]}

\bibitem[JK03]{JK03} N.~Joglar-Prieto, J. Koll\'ar,
 Real abelian varieties with many line bundles,
 Bull. London Math. Soc.~\textbf{35} (2003), 79--84

\bibitem[JM04]{JM04} N.~Joglar-Prieto, F.~Mangolte,
 Real algebraic morphisms and Del Pezzo surfaces of
   degree 2, J. Algebraic Geometry~\textbf{13}
 (2004), 269--285

\bibitem[Ko99]{KoIII} J.~Koll\'ar, Real algebraic
   threefolds III. Conic bundles,  J. Math. Sci., New
   York~\textbf{94} (1999), 996--1020

\bibitem[Ko00]{Ko00}
J.~Koll\'ar, 
Real algebraic
   threefolds IV. Del Pezzo fibrations, In: Complex analysis and
 algebraic geometry, de Gruyter, Berlin (2000), 317--346

\bibitem[KM08]{km} J.~Koll\'ar, F.~Mangolte,
Cremona transformations and homeomorphisms of  surfaces,
\texttt{arXiv:0809.3720 [math.AG]}

\bibitem[Ku99]{Ku} W.~Kucharz, 
Algebraic morphisms into
   rational real algebraic surfaces,
   J. Algebraic
   Geometry~\textbf{8} (1999), 569--579

\bibitem[Ma06]{Ma06} F.~Mangolte, 
Real algebraic
   morphisms on 2-dimensional conic bundles, Adv.\
   Geom.~\textbf{6} (2006), 199--213

\bibitem[Mu]{Mumford} D.~Mumford, Algebraic geometry I, Complex
 projective varieties (corrected second printing). Springer Verlag, 1976

\end{thebibliography}
\end{document}